\documentclass[12pt, authoryear]{amsart}
\usepackage[english]{babel}
\usepackage{amssymb,amsfonts,amsmath}
\usepackage{mathrsfs}
\usepackage{anysize}
\usepackage{url}
\usepackage[leqno]{amsmath}
\usepackage{tikz-cd}
\usepackage{stackengine}
\usetikzlibrary{decorations.pathreplacing}
 \usepackage[round]{natbib}

\usepackage{tikz-cd}
\usetikzlibrary{arrows}
\usetikzlibrary{patterns}
\usepackage{bbding}
\usepackage{datetime}
\usepackage{pgfplots}

\usepackage{caption}
\usepackage{tikz}
\usepackage{accents}
\usetikzlibrary{calc}
\usetikzlibrary{arrows}
\usetikzlibrary{backgrounds}

\newtheorem{theorem}{Theorem}[section]

\theoremstyle{definition}
\newtheorem{definition}[theorem]{Definition}

\usetikzlibrary{decorations.markings}
\tikzset{negated/.style={
        decoration={markings,
            mark= at position 0.5 with {
                \node[transform shape] (tempnode) {${\scriptstyle\setminus} $};
            }
        },
        postaction={decorate}
    }
}

\tikzset{degil/.style={
            decoration={markings,
            mark= at position 0.5 with {
                  \node[transform shape] (tempnode) {$\backslash$};
                  }
              },
              postaction={decorate}
}
}

\title[A chaotic Hopfield network]{A chaotic discrete-time continuous-state Hopfield  network with piecewise-affine activation functions}
\author{Benito Pires}
\email{benito@usp.br}

\address{Departamento de Computa\c c\~ao e Matem\'atica, Faculdade de Filosofia,
	Ci\^encias e Letras, \\ Universidade de S\~ao Paulo, Ribeir\~ao Preto, SP,
	14040-901, Brazil\\ benito@usp.br}

\begin{document}




\marginsize{2.5cm}{2.5cm}{1cm}{2cm}

\begin{abstract} We construct a chaotic discrete-time continuous-state Hopfield network 
with  piece\-wise-affine nonnegative activation functions and weight matrix with small positive entries. More precisely, there exists a Cantor set $C$ in the state space such that the network has sensitive dependence on initial conditions at initial states in $C$ and the network orbit of each initial state in $C$ has $C$ as its $\omega$-limit set. The approach we use is based on tools developed and employed recently in the study of the topological dynamics of piecewise-contractions. The parameters of the chaotic network are explicitly given. 


 \end{abstract}

  

\maketitle

\section{Introduction}

We consider nonnegative artificial neural networks (ANNs) consisting  of $n$ local units called \textit{neurons}, each of which takes as input a $n$-dimensional vector with nonnegative entries and uses a nonnegative activation function to process the weighted sum of the entries and  generate the output. The weights are the entries of a nonnegative matrix called \textit{weight matrix}. Nonnegative ANNs is an active field of research, see \cite{ Ali20171626,Ayinde20183969,Chorowski201462,Hosseini-Asl20162486}, \linebreak \cite{Lemme2012194,Su20181427}.

To understand the global dynamics of nonnegative ANNs, it is necessary to study the presence of attractors (e.g., fixed-points, cycles, fractal sets). Under the hypotheses that the activation functions are continuous and the state space is compact and convex, it follows from Brouwer's Fixed-Point Theorem that the network has at least one fixed-point. More generally, the existence of fixed-points in nonnegative ANNs with continuous activation functions was investigated in \cite{https://doi.org/10.48550/arxiv.2106.16239} by applying non-linear Perron-Frobenius theory.

In this article, we are concerned with nonnegative ANNs with discontinuous piecewise-affine activation functions. As we show here, discontinuities may result in  the existence of  fractal attractors. Since activation functions are supposed to map the entire real line into a small neighbourhood of two values (on/off), it is natural to assume that the activation functions are piecewise contractions.

The topological dynamics of piecewise contractions is an active field of research in dynamical systems theory, see, for instance, \cite{MR4120256,MR3394114,MR3225875,MR3820005,MR4030596}. By carefully choosing the activation functions, it is possible to embed the dynamics of  piecewise contractions of the interval into the network dynamics along an invariant line. With that approach, we can build nonnegative neural networks with a prescribed dynamics.         

The artificial neural networks considered in this article are known as Hopfield networks and  have been intensely studied since the pioneer works by   \cite{MR652033, WOS:A1984SU47600030}. To be more precise, a discrete-time continuous-state  Hopfield network (DCHN) consists of a  compact state space $X\subset [0,\infty)^n$ with non-empty interior and a piecewise-smooth map \linebreak $\mathbf{F}=(F_1,\ldots,F_n):[0,\infty)^n\to [0,\infty)^n$  such that $\mathbf{F}(X)\subset X$ and
 \begin{equation}\label{themapF}
F_i(x_1,\ldots,x_n)= f_i\left( \sum_{j=1}^n w_{ij } x_j-b_i\right)\quad (1\le i\le n),
\end{equation}
where $\mathbf{W}=(w_{ij})$ is a square matrix of size $n$ called \textit{weight matrix}, $\mathbf{b}=(b_1,\ldots,b_n)$ is a vector called \textit{external bias vector},
 and $f_1,\ldots,f_n$ are piecewise-smooth functions called \textit{activation functions}. Given an initial state $\mathbf{x}^{(0)}\in X$, the \textit{network state} at the time $k$ is the vector $\mathbf{x}^{(k)}\in X$ defined recursively by
\begin{equation}\label{xk} 
\mathbf{x}^{(k)}=\mathbf{F}\left(\mathbf{x}^{(k-1)}\right).
\end{equation}
In the terminology of neural networks, the update rule \eqref{xk} is called  \textit{syncronous} or \textit{parallel}.  
We call the whole sequence of network states
 $$(\mathbf{x}^{(k)})_{k\ge 0}=\left(\mathbf{x}^{(0)}, \mathbf{x}^{(1)},\mathbf{x}^{(2)},\ldots\right)$$ the \textit{orbit} of $\mathbf{x}^{(0)}$, which is completely determined by the initial state $\mathbf{x}^{(0)}$. Given an integer $p\ge 1$, we say that a sequence $\left(\mathbf{y}^{(k)}\right)_{k\ge 0}$ of $n$-dimensional vectors is a \textit{cycle of length $p$} if $\mathbf{y}^{(k+p)}=\mathbf{y}^{(k)}$  for all integers $k\ge 0$. A cycle $\left(\mathbf{x}^{(k)}\right)_{k\ge 0}$ such that $\mathbf{x}^{(k)}=\mathbf{F}\left(\mathbf{x}^{(k-1)}\right)$ for all $k\ge 0$  is called a \textit{network cycle}. Cycles of length $1$ are constant sequences. The orbit of $\mathbf{x}_0$ is a network cycle of length $1$ if and only if $\mathbf{x}^{(0)}$ is a fixed-point of the network, i.e., $\mathbf{F}\left(\mathbf{x}^{(0)}\right)=\mathbf{x}^{(0)}$. We say that the orbit of $\mathbf{x}_0$ is  \textit{asymptotic to a cycle} $\left(\mathbf{y}^{(k)}\right)_{k\ge 0}$ if
$\lim_{k\to\infty}\left\Vert \mathbf{x}^{(k)}-\mathbf{y}^{(k)}\right\Vert=0$, where $\Vert\cdot\Vert$ denotes the Euclidean norm in $\mathbb{R}^n$. We say that a Hopfield network is \textit{asymptotically periodic} if there is a finite collection of cycles such that each network orbit is asymptotic to a cycle of the collection.

It has been proved by  \cite{WOS:A1994PF90900007} that
 the existence of a Lyapunov function for a DCHN that decreases along the network orbits and is bounded from below plus the existence of an upper bound for the number of cycles imply that the network is asymptotically periodic.  He has also used a variant of the Lyapunov function provided in \cite{MR1102355} (see also \cite{WOS:A1989AE99600066}) to prove the following result.
   
 \begin{theorem}[\cite{WOS:A1994PF90900007}]\label{thm0} If a discrete-time continuous-state Hopfield network satisfies the hypotheses:
\begin{itemize}
\item [$(K1)$] The weight matrix $\mathbf{W}$ is symmetric and has non-negative diagonal;
\item [$(K2)$] The activation functions are increasing;
\item [$(K3)$] The number of network cycles is finite;
\end{itemize} 
then the network is asymptotically periodic and each network orbit is asymptotic to a cycle of length $1$ or $2$.
 \end{theorem}
 
    \cite{686695}  remarked that by combining \cite[Theorem 4]{WOS:A1994PF90900007} and the results in  \cite{21239}, it follows that the hypothesis $(K3)$ in Theorem \ref{thm0} is implied by $(K1)$ and $(K2)$ for an open dense set of pairs $(\mathbf{W},\mathbf{b})$ of symmetric weight matrices $\mathbf{W}$ and bias vectors $\mathbf{b}$.
  
In this article, we show that the hypothesis (K2) in Theorem \ref{thm0} is of paramount importance. More precisely, we provide an example of a chaotic discrete-time continuous-state Hopfield network that does not satisfy (K2).
 The approach we use is based on tools developed and employed recently in the study of the topological dynamics of piecewise contractions (see \cite{MR4120256,JGBP2022,MR3394114}).  The combination of discontinuity and contraction in the activation functions is the factor responsible for the rich dynamics consisting  of Cantor attractors in DCHNs.

\section{Statement of the result} 

We need some definitions from Chaos Theory to explain the dynamics of the Hopfield network. There are a variety of definitions of chaos (Devaney chaos, Li-Yorke chaos, {Wiggins chaos}, etc.).  The most important ingredient of chaos is the notion of sensitive dependence on initial conditions. 

\begin{definition}[sensitive dependence on initial conditions]\label{def0}
We say that the orbits 
 of a discrete-time continuous-state Hopfield network with state-space $X$ have  \textit{sensitive dependence on  initial conditions at $\mathbf{x}^{(0)}\in X$} if for some positive constant $\eta>0$ the following is true: for each $\epsilon>0$, there exist $\mathbf{y}^{(0)}\in X$ and $k\in\mathbb{N}$ such that $\left\Vert \mathbf{y}^{(0)}-\mathbf{x}^{(0)}\right\Vert\le\epsilon$ and
$\left\Vert \mathbf{y}^{(k)}-\mathbf{x}^{(k)}\right\Vert\ge\eta$, where $(\mathbf{y}^{(k)})_{k\ge 0}$ and $(\mathbf{x}^{(k)})_{k\ge 0}$ are the network orbits with initial states $\mathbf{y}^{(0)}$ and $\mathbf{x}^{(0)}$,  respectively.
\end{definition}

The \textit{$\omega$-limit set} of a network orbit $(\mathbf{x}^{(k)})_{k\ge 0}$ starting at the initial state $\mathbf{x}^{(0)}$ is the set
 \begin{equation}\label{limitset} 
 \omega(\mathbf{x}^{(0)})=
\omega\big((\mathbf{x}^{(k)})_{k\ge 0}\big)=\{\mathbf{p}\in \mathbb{R}^n: \exists  n_1<n_2<n_3<\cdots\,\,\textrm{such that}\,\,\lim_{k\to\infty}\mathbf{x}^{(n_k)}=\mathbf{p}\}.
 \end{equation}
  We are now ready to state the following definition.
 
\begin{definition}[chaotic Hopfield network]\label{cfp} We say that a discrete-time continuous-state Hopfield network with state space $X$ is \textit{chaotic} if there exists a Cantor set $C\subset X$ such that the following statements are true:
\begin{itemize}
\item [(C1)] $\omega\big(\mathbf{x}^{(0)}\big)=C$ for each $\mathbf{x}^{(0)}\in C$,
\item [(C2)] The network orbits have sensitive dependence on initial conditions at all $\mathbf{x}^{(0)}\in C$.
\end{itemize}
\end{definition}

Now we introduce the parameters of the chaotic Hopfield network. {In what follows, we denote by {$\boldsymbol{\omega}=\omega_0\omega_1\omega_2\ldots$} the \textit{Fibonacci word}, that is, the sequence of binary digits
\begin{equation}\label{fw}
 \boldsymbol{\omega}=010010100100101001010010010100100101001010010010100 . . .
\end{equation}
defined by {$\omega_i=2+\lfloor (i+1)\varphi\rfloor-\lfloor (i+2)\varphi\rfloor$}, where $\varphi=(1+\sqrt{5})/2$ is the golden ratio and $\lfloor x\rfloor$ denotes the integral part of $x$. Notice that the Fibonacci word is the sequence A003849 in the OEIS\footnote{\textsc{The on line encyclopedia of integer sequences}\textsuperscript{\textregistered}.}.

Now we are ready to state our main result.

\begin{theorem}[main result]\label{theorem1} Let $\mathscr{H}$ be a discrete-time continuous-state Hopfield network satisfying the following conditions:
\begin{itemize}
\item [$(H1)$] The state space is $X=[0,1]^n$ and the external bias satifies $b_i=0$ for each $i$;
\item [$(H2)$] The weight matrix $\mathbf{W}=(w_{ij})$ has positive entries and $\frac34<\sum_{j=1}^n w_{ij}<1$, $\forall i$;

\item [$(H3)$] The activation functions $f_i:[0,\infty)\to [0, \infty)$ are defined by   
$$ f_i(x) =  
 \begin{cases} \dfrac{1}{2\rho} x + \delta v_i   & \textrm{if}\quad x \in  \left[0, 2(1-\delta)\rho v_i\right), \\[0.15in] 
   \dfrac{1}{2\rho} x + (\delta-1)v_i  & \textrm{if} \quad x\in \left[2(1-\delta)\rho v_i, \rho v_i\right], \\[0.15in]
    {\dfrac{1}{2} }v_i + (\delta -1) v_i & \textrm{if} \quad x\in [\rho v_i,\infty),
 \end{cases}  \quad\textrm{where}\quad   \delta = \dfrac12 + \dfrac14\sum_{k\ge 0} \dfrac{\omega_k}{2^k},                
 $$
 $\boldsymbol{\omega}=\omega_0\omega_1\omega_2\ldots$ is the Fibonacci word defined in \eqref{fw}, $\rho$ is the Perron-Frobenius eigenvalue of $\mathbf{W}$, and $\mathbf{v}=(v_1,\ldots,v_n)$ is the associated probability eigenvector.
\end{itemize}
Then $\mathscr{H}$ is a chaotic Hopfield network.
  \end{theorem}
Below is an activation function satisfying the hypothesis (H3). Notice that the activation function is not increasing, therefore it does not satisfy (K2) in Theorem \ref{thm0}.
  
  \begin{figure}[!htb]
        \centering
\fbox{
    \begin{tikzpicture}[scale=0.4]



\draw [  thick] (-5,0) -- (2.45,0);
\draw [  thick, ->] (-5,0) -- (6,0) node [below] {  $x$};
\draw [  thick, ->] (-5.1,0) -- (-5.1,6) node [above] {  $f_i(x)$};
	
			

\draw[-, very thick, red] (-5,2.3)--(-1.5,4);
\draw[-, very thick, red] (-1.5,0.7)--(0,1.3);
\draw[-, very thick, red] (-0,1.3)--(5,1.3);

\draw[black] node at (2, 6) {$\theta=2(1-\delta)\rho$};



\draw  (-1.5,0) node [below, xshift=-0.2cm]{$\theta v_i$};
\draw  (0,0) node [below, xshift=0.1cm,  yshift=-0.1cm]{$\rho v_i$};
\draw[thick] (0,-0.2)--(0,0.2);
\draw[thick] (-1.5,-0.2)--(-1.5,0.2);

\draw[white] node at (-1.5,4) {\textbullet};
\draw[red] node at (-1.5,4) {$\circ$};
\draw[red] node at (-1.5,0.7) {\textbullet};

\end{tikzpicture}}

  \caption{Example of activation function in Theorem \ref{theorem1}}\label{fig3}
\end{figure}
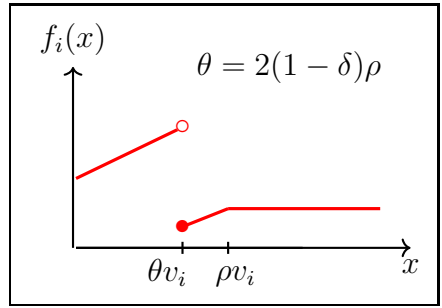

 \section{Proof of the main result}
 
  Let $\mathscr{H}$ be a discrete-time continuous-state Hopfield network satisfying conditions (H1)-(H3) in Theorem \ref{theorem1}. We keep all the notation introduced in the statement of Theorem \ref{theorem1}.  The network dynamics is ruled by $\mathbf{F}=(F_1,\ldots,F_n):[0,\infty)^n\to [0,\infty)^n$ defined by 
  \begin{equation}\label{themapF2}
F_i(x_1,\ldots,x_n)= f_i\left( \sum_{j=1}^n w_{ij } x_j\right)\quad (1\le i\le n).
\end{equation}
The proof follows from a sequence of claims. First we need to verify that the activation functions $f_i$ are well-defined 
and  $\mathbf{F}$ takes the state space $X= [0, 1]^n$ into itself.\\

\noindent Claim A. $\frac12< 2(1-\delta)<\rho <1$.\\

In fact, by the definition of $\delta$ and $\boldsymbol{\omega}$, we have that $\frac58<\delta<\frac34$, thus $\frac12<2(1-\delta)<\frac34$. Since $\rho$ is the spectral radius of $\mathbf{W}$, by (H2) and by Lemma 2.8 in \cite[p. 36]{MR1753713}, it follows that
$\frac 34 < \rho < 1$. In this way, $\frac12<2(1-\delta)  < \rho<1$. \\

\noindent Claim B. $F(X)\subset X$.\\

Given $\mathbf{x}=(x_1,\ldots, x_n)\in [0,1]^n$ and $1\le i\le n$, let $y_i = \sum_{j=1}^n w_{ij} x_j$. We claim that $F_i(\mathbf{x})\le v_i \le 1$. In fact, by (H3),  \eqref{themapF2} and Claim A, we have that
$$ F_i(\mathbf{x})= f_i (y_i) = \begin{cases} \dfrac{1}{2\rho} y_i + \delta v_i < \dfrac{2(1-\delta)\rho v_i}{2\rho} + \delta v_i = v_i& \text{if}\quad y_i\in \left[0, 2(1-\delta)\rho v_i \right) \\[0.1in]
 \dfrac{1}{2\rho} y_i + (\delta -1) v_i < \dfrac{\rho v_i}{2\rho} + (\delta-1) v_i< v_i& \text{if}\quad y_i\in \left[2(1-\delta)\rho v_i, \rho v_i \right] \\[0.1in]
 \dfrac{1}{2} v_i + (\delta-1) v_i < \dfrac12 v_i -\dfrac14 v_i < v_i & \text{if}\quad y_i\in \left[\rho v_i, \infty \right)
\end{cases}.
$$
Likewise, we claim that $F_i(\mathbf{x})>0$. In fact,
$$ F_i(\mathbf{x})= f_i (y_i) = \begin{cases} \dfrac{1}{2\rho} y_i + \delta v_i \ge 0& \text{if}\quad y_i\in \left[0, 2(1-\delta)\rho v_i \right) \\[0.1in]
 \dfrac{1}{2\rho} y_i + (\delta -1) v_i \ge \dfrac{2(1-\delta) \rho v_i}{2\rho} + (\delta-1) v_i= 0& \text{if}\quad y_i\in \left[2(1-\delta)\rho v_i, \rho v_i \right] \\[0.1in]
 \dfrac{1}{2} v_i + (\delta-1) v_i \ge  \dfrac12 v_i -\dfrac12 v_i \ge 0 & \text{if}\quad y_i\in \left[\rho v_i, \infty\right)
\end{cases}.
$$
We have proved that ${F}_i\big([0,\infty)^n\big)\subset [0,v_i]\subset [0,1]$ for each $1\le i\le n$. In particular, we have that
$\mathbf{F}$ takes the state space $X=[0,1]^n$ into itself.\\

To show that $\mathscr{H}$ is chaotic, it suffices to study the dynamics of $\mathbf{F}$ along the half-line
$$ L= \{t\mathbf{v}: t\ge 0\}.
$$  

\noindent Claim C. $L$ is $\mathbf{F}$-invariant, that is, $\mathbf{F}(L)\subset L$.\\

 In fact, by \eqref{themapF} and  (H3), if follows that
 if $\{\mathbf{e}_1,\ldots,\mathbf{e}_n\}$  denotes the canonical basis of $\mathbb{R}^n$, then for each
$\mathbf{x}=t\mathbf{v}\in L$, we have that
\begin{equation}\label{Fx}
 \mathbf{F}(\mathbf{x})=\sum_{i=1}^n  f_i\left( \sum_{j=1}^n w_{ij } t v_j\right)\mathbf{e}_i=\sum_{i=1}^n f_i(\rho t v_i)\mathbf{e}_i=
 \sum_{i=1}^n  v_i \widetilde{g}(t) \mathbf{e}_i = \widetilde{g}(t) \mathbf{} \mathbf{v},
\end{equation}
where $\widetilde{g}:[0,\infty) \to [0,1]$ is the piecewise-affine map defined by
\begin{equation}\label{gg} 
\widetilde{g}(t)  = \begin{cases} \frac12 t + \delta & \textrm{if}\quad t \in \big[0, 2(1-\delta)),\\[0.1in]
\frac12 t + \delta - 1 & \textrm{if}\quad t \in \big[ 2(1-\delta), 1],\\[0.1in]
 \delta - \frac12& \textrm{if}\quad t \in \big[ 1, \infty).
\end{cases}
\end{equation}
In this way, $\mathbf{F}(L)\subset L$, which proves the claim.\\

In what follows, let $h:L\to [0,\infty)$ be the homeomorphism defined by $h(s\mathbf{v})=s$, $s\ge 0$, and $\widetilde{g} :[0,\infty)\to [0,1]$ be as in \eqref{gg}. \\

\noindent Claim D. The following diagram commutes
\[ \begin{tikzcd}[arrows={-Stealth}]
{L}\rar["\mathbf{F}\vert_L"]\dar["h"] & {L}\dar["h"]   \\%
{[0,\infty)}\rar[swap, "\widetilde{g}"] & {[0,\infty)}  
\end{tikzcd},
\]
that is,  $h\big(\mathbf{F}(\mathbf{x})\big)=\widetilde{g}\big(h(\mathbf{x})\big)$ for all $\mathbf{x}\in L$.\\

In fact, by \eqref{Fx}, for all $\mathbf{x}=t \mathbf{v}\in L$, we have that
$$
h\big(\mathbf{F}(\mathbf{x})\big)= h\big(\widetilde{g}(t)\mathbf{v}\big)= \widetilde{g}(t) = \widetilde{g}\big(h(t\mathbf{v})\big)= \widetilde{g}\big(h(\mathbf{x})\big).
$$

In the terminology of dynamical systems, Claim D states that the restriction  $\mathbf{F}\vert_{L}$ to the invariant set $L$
 and the piecewise-contraction $\widetilde{g}:[0,\infty)\to [0,1]$ defined in \eqref{gg} are topologically conjugate by the conjugacy $h$. In this way, the maps $\mathbf{F}\vert_L$ and $\widetilde{g}$ have the same topological dynamics. The dynamics of the map $\widetilde{g}:[0,\infty)\to [0,1]$ is completely known. More precisely, $\widetilde{g}\big([0,\infty)\big)\subset [0,1]$ and 
 the restriction $g=\widetilde{g}\vert_{[0,1]}:[0,1]\to [0,1]$ is the piecewise-contraction with  one discontinuity defined by
 \begin{equation}\label{ggg} 
{g}(t)  = \begin{cases} \frac12 t + \delta & \textrm{if}\quad t \in \big[0, 2(1-\delta)),\\[0.1in]
\frac12 t + \delta - 1 & \textrm{if}\quad t \in \big[ 2(1-\delta), 1].
\end{cases}
\end{equation}
The dynamics of $g$ was studied in \cite{JGBP2022}. To conclude the proof, we will verify that there exists a Cantor set $C_L\subset L$ such that
 the Hopfield network $\mathscr{H}$ satisfies Conditions (C1) and (C2) in Definition \ref{cfp}.\\

\noindent Claim E. The Hopfield network $\mathscr{H}$ satisfies Conditions (C1) in Definition \ref{cfp}.\\

In fact, by \cite[Lemma 6]{JGBP2022} with $b=2$, there exists a Cantor set $C\subset [0,1]$ such that 
$\omega_g\big(t^{(0)}\big)=C$ for all $t^{(0)}\in C$, where $\omega_g\big(t^{(0)}\big)$ denotes the $\omega$-limit set of $t^{(0)}$ by the map $g$ defined in \eqref{ggg}. More precisely, 
$$ \omega_g\big(t^{(0)}\big)=\{p\in [0,1]: \exists  n_1<n_2<n_3<\cdots\,\,\textrm{such that}\,\,\lim_{k\to\infty}t^{(n_k)}=p\},$$
where $\big(t^{(k)}\big)_{k\ge 0}$ is the $g$-orbit of $t_0$ defined recursively by $t^{(k)}=g\big(t^{(k-1)}\big)$. Since $\widetilde{g}\vert_{C}=g\vert_{C}$, where $\widetilde{g}$ is as in \eqref{gg}, we have that $\omega_{\widetilde{g}}\big(t^{(0)}\big)=C$ for all $t^{(0)}\in C$. Now let $C_L\subset L$ be the Cantor set defined by $C_L=h^{-1}(C)$, where $h$ is the homeomorphism in Claim D. Let $\mathbf{x}^{(0)}\in C_L$, then there exists a unique $t_0 \in C$ such that $\mathbf{x}^{(0)}= h^{-1}\big(t^{(0)}\big)$. By Claim D, it follows that
$\omega\big(\mathbf{x}^{(0)}\big)=C_L$, where $\omega\big(\mathbf{x}^{(0)}\big)$ denotes the $\omega$-limit set of $\mathbf{x}^{(0)}$ defined in \eqref{limitset}. In this way, Condition (C1) in Definition \ref{cfp} holds true.\\

\noindent Claim F. The Hopfield network $\mathscr{H}$ satisfies Conditions (C2) in Definition \ref{cfp}.\\

We will use the same notation introduced in the proof of Claim E. By Claim D, the Hopfield network $\mathscr{H}$ has sensitive dependence on initial conditions at a point
$\mathbf{x}^{(0)}=h^{-1}\big(t^{(0)}\big)$ of the Cantor set $C_L$ if and only if the interval map $g:[0,1]\to [0,1]$ has the same property at the point $t^{(0)}$ of the Cantor set $C$, that is, if for some positive constant $\eta>0$, the following is true: for each $\epsilon>0$, there exists $s^{(0)}\in [0,1]$ and $k\in\mathbb{N}$ such that $\left\vert s^{(0)} - t^{(0)}\right\vert\le\epsilon$ and $\left\vert s^{(k)} - t^{(k)}\right\vert \ge \eta$, where $\big(s^{(k)}\big)_{k\ge 0}$ and
$\big(t^{(k)}\big)_{k\ge 0}$ are the $g$-orbits of $s^{(0)}$ and $t^{(0)}$, respectively. 

To conclude the proof of Claim F, let us show that $g$ has sensitive dependence on initial conditions at any $t^{(0)}\in C$. 
Without loss of generality, we may assume that $t^{(0)}>0$. Given $\epsilon>0$, let $0<\epsilon'<\epsilon$ be so small that
 $J=(t_0-\epsilon', t_0+\epsilon')$ is an open subinterval of  $[0,1]$. We have two cases to consider.\\
 
 \noindent Case I. The discontinuity $d=2(1-\delta)$ belongs to $g^k(J)$ for some $k\ge 0$. \\
 
 In this case, if $k$ is the least nonnegative integer such that $d\in g^k(J)$, then $g^k(J)$ is an open interval containing $d$. Moreover, there exist $s^{(0)}\in J$ such that one of the following alternatives occurs:
 \begin{itemize}
 \item [$(i)$]
$g^k\big(s^{(0)}\big)<d$ and  $g^k\big(t^{(0)}\big)\ge d$;
\item [$(ii)$]
$g^k\big(s^{(0)}\big)>d$ and  $g^k\big(t^{(0)}\big)<d$.
\end{itemize}
 It is elementary to verify that in either case $\left\vert g^{k+1}\big(s^{(0)}\big)-g^{k+1}\big(t^{(0)}\big)\right\vert\ge \frac12$. Moreover, since $s^{(0)}\in (t_0-\epsilon',t_0+\epsilon')$, we have that $\left\vert s^{(0)}- t^{(0)}\right\vert\le \epsilon'<\epsilon$. In this way, $g$ has sensitive dependence on initial conditions at $t^{(0)}\in C$. By the previous discussion, the network $\mathscr{H}$ has sensitive dependence on initial conditions at the point $\mathbf{x}^{(0)}=h^{-1}\big(t^{(0)}\big)\in C_L$. Since $t^{(0)}$ is an arbitrary point of ${C}$, we have that $\mathbf{x}^{(0)}$ is an arbitrary point of $C_L$. This concludes the proof of Case I.
  \\ 
 
  \noindent Case II. The discontinuity $d=2(1-\delta)$ does not belong to $\bigcup_{k\ge 0} g^k(J)$.\\ 
  
  Let $\alpha = (3 - \sqrt{5})/2$. By \cite[Lemma 6]{JGBP2022}, there exists a continuous, nondecreasing and surjective map  $u:[0,1]\to [0,1]$ (called topological conjugacy) such that 
  \begin{equation}\label{ut}
 u(t)=1-\alpha \,\,\,\, \textrm{if and only if} \,\,\,\, t = d
  \end{equation}
   and
    the following diagram commutes 
\begin{equation}\label{dc00}
\begin{tikzcd}[arrows={-Stealth}]
{[0,1]}\rar["g"]\dar["u"] & {[0,1]}\dar["u"]   \\%
{[0,1]}\rar[swap, "T"] & {[0,1]}  
\end{tikzcd},
\end{equation}
that is, $u\circ g=T\circ u$, where $T:[0,1]\to [0,1]$ is the interval map (called minimal interval exchange transformation or irrational rotation by $\alpha$) defined by
$$T(t) = \begin{cases} t + \alpha & \textrm{if}\quad t\in [0, 1-\alpha)\\
t + \alpha - 1 & \textrm{if} \quad t \in [1 - \alpha, 1]
\end{cases}. 
$$
Since $T$ is equivalent to the irrational rotation by  $\alpha$, we have that $T^{-1}$ is equivalent to the irrational rotation by  $-\alpha$. In particular, every orbit of $T^{-1}$ is dense. Hence, given any interval $U\subset [0,1]$ of positive length, there exists an inteter
$k\ge 0$ such that $T^{-k}(1-\alpha)\in U$, or equivalently, $1-\alpha\in T^k(U)$. Moreover, by the item $(iii)$ of \cite[Lemma 6]{JGBP2022}, it follows that $u\big({J}\big)$ is  an interval of positive length. In this way, by all the previous discussion and by \eqref{dc00}, there exists an integer $k\ge 0$ such that $1-\alpha\in T^k\big( u(J)\big)=u\big(g^k(J)\big)$. By \eqref{ut}, $d\in g^k(J)$, which contradicts the hypothesis of Case II. Hence, Case II cannot occur.

  \section{Concluding Remarks}
  
  In this article, we have used  modern techniques 
   from the area of Dynamical Systems to construct a chaotic discrete-time continuous-state Hopfield network whose parameters are given explicitly. The activation functions used are non-increasing piecewise affine-contractions. Chaotic Hopfield networks are rare and difficult to construct because some of their parameters in general are transcendental numbers. In this article, readers will find a formal mathematical proof that the Hopfield network built here is chaotic. 
 
We have considered the set $X=[0,1]^n$ as the state space of the Hopfield network $\mathscr{H}$ in the sense that the map $\mathbf{F}$ defined by \eqref{themapF} takes $X$ into $X$. The  proof that we provide shows that, in fact, $\mathbf{F}$ takes the $n$-dimensional rectangle $[0,v_1] \times \cdots\times [0,v_n]$ into itself, so that we can also consider $X=\Pi_{i=1}^n [0,v_i]$ as the state space, where $(v_1,\ldots, v_n)$ is the Perron-Frobenius probability eigenvector of  the weight matrix $\mathbf{W}$. In either case, the state space is compact and therefore all the orbits are bounded. Moreover, all the network orbits with the initial state $\mathbf{x}^{(0)}$ in some Cantor set $C$ has the Cantor set $C$ as its $\omega$-limit set. In this case, it is not possible to predict the network state $\mathbf{x}^{(k)}$ for $k$ large because of the sensitive dependence on initial conditions.

The hypothesis (H2) in Theorem \ref{theorem1} that the  $i$-th row of the weight matrix $\mathbf{W}=(w_{ij})$ satisfies $\frac34<\sum_{j=1}^n w_{ij}<1$ for all $i$, can be replaced by the hypothesis that the $j$-th column of the weight matrix satisfies $\frac34<\sum_{i=1}^n w_{ij}<1$ for all $j$.

 \section*{Declaration of competing interest}
 
 The author declares that he has no known competing financial interests or personal relationships that could have appeared to influence the work reported in this paper.

 \section*{Acknowledgments.} 
The author was partially supported by grant {\#}2019/10269-3 S\~ao Paulo Research Foundation (FAPESP). The author is very thankful to Antonio C. Roque for fruitful discussions
 on neural networks. The author is also grateful to Aluizio F. R. Ara\'ujo for introducing him to  Hopfield Networks  twenty years ago.


\end{document}